\documentclass[11pt]{amsart}
\usepackage{fullpage}
\usepackage{latexsym}
\usepackage{amssymb}

\input amssym.def
\input amssym.tex

\def\page#1{\leaders\hbox to 5mm{.}\hfill \rlap{\hbox to 2mm{\hfill #1}}\par}

\font\tengoth=eufm10 at 12pt
\font\sevengoth=eufm7
\font\fivegoth=eufm5
\newfam\gothfam
\textfont\gothfam=\tengoth
\scriptfont\gothfam=\sevengoth
\scriptscriptfont\gothfam=\fivegoth

\def\thepec{2mmplus1mmminus1mm} 
\def\thegec{4mmplus1mmminus1mm} 
\def\pec{\vskip\thepec}
\def\gec{\vskip\thegec}

\def\og{\leavevmode\raise.3ex\hbox{$\scriptscriptstyle\langle\!\langle\,$}}
\def\fg{\leavevmode\raise.3ex\hbox{$\scriptscriptstyle\,\rangle\!\rangle\ $}}

\newtheorem{theorem}{Theorem}
\newtheorem{cor}[theorem]{Corollary}
\newtheorem{lemma}[theorem]{Lemma}
\newtheorem{prop}[theorem]{Proposition}
\newtheorem{quest}[theorem]{Question}
\newtheorem{conj}[theorem]{Conjecture}

\newtheorem{pb}[theorem]{Problem}
\newtheorem{fact}[theorem]{Fact}

\begin{document}

\title{Kac-Moody groups as discrete groups}
\author{Bertrand R\'emy}
\maketitle

\pec 

{\footnotesize
{\sc Abstract}.---~
This survey paper presents the discrete group viewpoint on Kac-Moody groups. 
Over finite fields, the latter groups are finitely generated; they act on new
buildings enjoying remarkable negative curvature properties. 
The study of these groups is shared between proving results supporting the analogy with
some $S$-arithmetic groups, and exhibiting properties showing that they are new groups.}

\pec

{\footnotesize

\pec 

}

\vskip 10mm 

\section*{Introduction}

Kac-Moody groups were initially designed to generalize algebraic groups \cite{Tit87}. 
They share many properties with the latter groups, mainly of combinatorial nature. 
For instance, they admit $BN$-pair structures \cite[IV.2]{Bou81}, and actually a much
finer combinatorial structure -- called {\it twin root datum~} -- formalizing the existence
of root subgroups permuted by a (possibly infinite) Coxeter group \cite{Tit92} (survey papers
on this are for instance \cite{Tit89} and \cite{RemBie}). 
Moreover Kac-Moody groups and some twisted versions of them are expected to be the group side
in the classification of a reasonable class of buildings, the so-called Moufang 2-spherical
twin buildings. 
Our intention is not to go into detail about this; nevertheless, we shortly recalled these
facts in order to emphasize that the goal of this paper is to present a significant change of
viewpoint on Kac-Moody groups. 

\pec 

The geometric counterpart to the group combinatorics of a $BN$-pair is the existence
of an action on a remarkable geometry: a building. 
To any Coxeter system $(W,S)$ is attached a simplicial complex $\Sigma$ on the
maximal simplexes of which the group $W$ acts simply transitively \cite[\S 2]{RonLec}. 
A {\it building~} is a simplicial complex, covered by subcomplexes all isomorphic to the
same $\Sigma$ -- called {\it apartments}, and satisfying remarkable incidence
properties:  any two {\it facets~} (i.e. simplices) are contained in an apartment, and given
any two apartments $A, A' \simeq \Sigma$, there is a simplicial isomorphism $A \simeq A'$
fixing $A \cap A'$ -- see \cite[p. 77]{BroBuildings}, and also \cite[\S 3]{RonLec} for the
chamber system approach. 

\pec 

In order to define a general Kac-Moody group $\Lambda$, we need a integral matrix $A$
satisfying properties much weaker than those defining Cartan matrices in the classical sense
of complex semisimple Lie algebras \cite[Introduction]{Tit87}. 
We also need to choose a ground field (which will always be a finite field ${\bf F}_q$ in
what follows). 
Writing down a presentation of $\Lambda$ would require a lot of combinatorial and Lie
algebra material we wouldn't use later -- see \cite[Subsect. 3.6]{Tit87} and 
\cite[Sect. 9]{RemAst} for details. 
It is a basic result of the theory that to a Kac-Moody group is naturally attached a pair of
twin buildings (via $BN$-pairs) on the product of which it acts diagonally 
(Fact \ref{fact - buildings} of the present paper). 
The standard example of such a group is $\Lambda={\bf G}\bigl( {\bf K}[t,t^{-1}] \bigr)$ for
${\bf G}$ a semisimple group over a field ${\bf K}$. 
Now, if we choose the ground field to be ${\bf F}_q$, the latter example is
an arithmetic group in positive characteristic, and in this case the buildings alluded to
above are the Bruhat-Tits buildings of the non Archimedean semisimple Lie groups 
${\bf G}\bigl({\bf F}_q(\!(t)\!)\bigr)$ and ${\bf G}\bigl({\bf F}_q(\!(t^{-1})\!)\bigr)$. 

\pec 

The latter examples are a very special case of Kac-Moody groups, called of {\it affine
type}, but many other cases are available. 
On the one hand, this suggests to see Kac-Moody groups as generalized arithmetic groups
over function fields, which leads to natural questions, e.g. asking whether some classical
properties of discrete subgroups of Lie groups are relavant or true. 
On the other hand, it can be shown that some new buildings can be produced thanks to
Kac-Moody groups, and this leads to asking whether the groups attached to exotic buildings
are themselves new. 
Note that \og new\fg in this context means that the buildings are neither of spherical nor
of affine type, i.e. don't come from the classical Borel-Tits (resp. Bruhat-Tits) theory on
algebraic groups over arbitrary (resp. local) fields. 

\pec 

From the point of view of metric spaces, Kac-Moody theory is an algebraic way to construct
spaces with non-positively curved, often hyperbolic, distances and admitting highly
transitive isometry groups. 
The algebraic origin of these groups enables to obtain interesting structure results for
various isometry groups (discrete or much bigger). 
Therefore, in the case of hyperbolic Kac-Moody buildings, techniques from group
combinatorics such as $BN$-pairs and from hyperbolic spaces \og \`a la Gromov\fg can be
combined. 
This leads us to say that the general trend to study finitely generated Kac-Moody groups is
from algebraic and combinatorial methods to geometric and dynamical ones. 
In this paper, we explain for instance how the theory of finitely generated Kac-Moody
groups, i.e. Kac-Moody groups over finite fields, naturally leads to studying uncountable
totally disconnected groups generalizing semsimple groups over local fields of positive
characteristic, groups which we call {\it topological Kac-Moody groups}. 
We are especially interested in proving that the groups we obtain are new in general, by
proving that they cannot be linear over any field. 
In short, the study of finitely generated Kac-Moody groups is shared between proving
classical properties by comparing them to (linear) lattices of non Archimedean Lie groups,
and finally standing by a difference with the classical situation to disprove linearities. 

\pec 

This appraoch is not new, since for lattices of products of trees M. Burger and Sh. Mozes
managed to prove many classical properties of lattices (among which the normal subgroup
theorem) but proved that an important difference is the possibility to obtain non residually
finite groups, which is impossible for finitely generated linear groups \cite{BurMozCras},
\cite{BurMozTrees}, \cite{BurMozProd}. 
The main application is the construction of the first finitely generated torsion free simple 
groups. 
We note also that Y. Shalom's work \cite{Shalom} shows that thanks to representation
theory, many properties can be proved for irreducible lattices of products of general
locally compact groups. 
Since trees are one-dimensional special cases of buildings, the previous
references encourage us to think that finitely generated Kac-Moody groups will produce
interesting examples of groups which are not linear but handable via their diagonal
actions on products of buildings. 

\pec 

This paper is organized as follows. 
In Sect. \ref{s - arithmetic}, it is shown why Kac-Moody groups over finite fields should
be seen as generalizations of some $S$-arithmetic groups in positive characteristic. 
It is also explained how they provide new buildings and why Kac-Moody groups should be
expected to be new groups. 
In Sect. \ref{s - generalized KM}, we are interested in a specific class of hyperbolic
buildings. 
In this context, we can produce non-isomorphic Kac-Moody groups with the same buildings, and 
discrete groups which are close to Kac-Moody groups, but with several ground fields: they
have strong non-linearity properties.  
In Sect. \ref{s - topological}, we are interested in totally disconnected groups generalizing
semisimple groups over local fields, arising as closures of non-discrete actions of Kac-Moody
groups on buildings.  
We quote the existence of a nice combinatorial structure, as well as a
topological simplicity result for them. 
In Sect. \ref{s - hard NL}, we sketch the proof of complete non-linearity of some Kac-Moody 
groups. 
This is where we use the topological groups of the previous section. 
We actually mention that there are some Kac-Moody groups all of whose linear images are
finite, whatever the target field. 
In Sect. \ref{s - conjectures}, we ask some questions about the various groups previously
defined in the paper. 
We conjecture the non-linearity of a wide class of finitely generated Kac-Moody groups, and
the non-amenablity as well as the abstract simplicity of topological Kac-Moody groups. 

\pec 

The author expresses his deep gratitude to the organizers of the conference \og
Geometric Group Theory\fg (Guwahati, Assam -- India), organized by the Indian Institute of
Technology Guwahati and the Indian Statistical Institute, and supported by the National Board
of Higher Mathematics. 
Meenaxi Bhattacharjee managed to make a tremendous human and scientific event out of a
mathematical conference. 

\gec 

\section{Generalized arithmetic groups acting on new buildings}
\label{s - arithmetic}

This section is mainly dedicated to quoting results supporting the analogy between Kac-Moody
groups over finite fields and $\{ 0;\infty \}$-arithmetic groups over function fields 
\cite[Sect. 2-3]{RemNewton}.  
The arguments are: the existence of a discrete diagonal action on a product
of two buildings (\ref{ss - lattice}), cohomological finiteness properties 
(\ref{ss - finiteness}) and continuous cohomology vanishings  (\ref{ss - T}).  
In the last two subsections  (\ref{ss - hyperbolic} and \ref{ss  - easy NL}), we provide
arguments showing that Kac-Moody theory does provide interesting new group-theoretic/geometric
situations. 

\subsection{Discrete actions on buildings and finite covolume} 
\label{ss - lattice}
A Kac-Moody group is generated by those of its root groups which are indexed by simple
roots and their opposites (in finite number), and by a suitable maximal torus normalizing
them \cite[\S 3.6]{Tit87}. 
Since over ${\bf F}_q$ all these groups are finite, we obtain: 

\begin{fact}
\label{fact - finitely generated}
Any Kac-Moody group $\Lambda$ over any finite field is finitely generated. 
\end{fact}

From now on, $\Lambda$ is a Kac-Moody group defined over ${\bf F}_q$. 
Recall that a group action on a building is {\it strongly transitive~}if it is transitive on 
the inclusions of a chamber in an apartment. 
Combining \cite[Subsect. 5.8, Proposition 4]{Tit87} and \cite[Theorem 5.2]{RonLec}, we
obtain: 

\begin{fact}
\label{fact - buildings}
To $\Lambda$ are attached two isomorphic, locally finite buildings $X_\pm$, each of them
admitting a strongly transitive $\Lambda$-action. 
\end{fact}

This fact is fundamental to understand Kac-Moody groups: the geometry of the buildings is the
basic substitute for a natural structure on 
$\Lambda$ arising from infinite-dimensional algebraic geometry. 
In the specific case of an $S$-arithmetic group $\Lambda={\bf G}\bigl( {\bf F}_q[t,t^{-1}]
\bigr)$ (for $S=\{0;\infty\}$), the buildings $X_\pm$ are the Bruhat-Tits buildings of ${\bf
G}\bigl( {\bf F}_q (\!( t )\!) \bigr)$ and ${\bf G}\bigl( {\bf F}_q (\!( t^{-1} )\!) \bigr)$. 
In general, it is still true that for an arbitrary Kac-Moody group, the diagonal action on
the product of buildings $X_- \times X_+$ is discrete \cite{RemCras}. 
Moreover the action has a nice convex fundamental domain contained in a single apartment and
defined as an intersection of roots (seen as half-apartments) 
\cite[\S 3, Corollary 1]{Abr97}. 
The next step in the analogy with arithmetic groups consists in asking whether the finitely
generated group $\Lambda$ is a {\it lattice~} of $X \times X_-$, meaning that the locally
compact group ${\rm Aut}(X_-) \times {\rm Aut}(X_+)$ moded out by the image of
$\Lambda$ carries a finite invariant measure \cite[0.40]{Margulis}. 
This is the main result of \cite{RemCras}: 

\begin{theorem}
\label{th - lattice}
Assume the Weyl group $W$ of $\Lambda$ is infinite and  denote by 
$W(t):=\sum_{w \in W}t^{\ell(w)}$ its growth series. 
Assume that $W({1\over q}) < \infty$.
Then $\Lambda$ is a lattice of $X \times X_-$ for its diagonal action, 
and for any point $x_- \! \in \! X_-$ the stabilizer $\Lambda(x_-)$ is a lattice of $X$.
These lattices are never cocompact.
\end{theorem}

In the arithmetic case of  ${\bf G}(\mathbf{F}_q[t,t^{-1}])$, the Weyl group has polynomial
growth (it is virtually abelian because it is an affine reflection group), so the condition 
$W({1 \over q}) < \infty$ is empty and ${\bf G}(\mathbf{F}_q[t,t^{-1}])$ is always a
lattice of ${\bf G}\bigl( {\bf F}_q (\!( t )\!) \bigr) \times 
{\bf G}\bigl( {\bf F}_q (\!( t^{-1} )\!) \bigr)$, whatever the value of the prime power $q$. 
This is a particular case of a well-known result in reduction theory in positive
characteristic \cite{BehrRed}, \cite{Harder}. 

\subsection{Finiteness properties} 
\label{ss - finiteness}
Cohomological finiteness properties is a very hard problem for arithmetic groups in positive
characteristic \cite{Behr}, which makes a sharp difference with the number field case. 
Nevertheless, it is natural to expect that some results, similar to those which are known
for arithmetic groups in the function field case, should hold for Kac-Moody groups. 
This is indeed the case, up to taking into account more carefully the submatrices in the 
generalized Cartan matrix defining the group -- this is the theme of P.
Abramenko's book \cite{Abr97}. 
The following theorem sums up Theorems 1 and 2 from \cite{AbrBie} in a slightly
different language. 
For instance, we introduce the parabolic subgroups (defined by means of $BN$-pairs in 
[loc. cit.]) via the group action on the buildings: such a subgroup is a facet fixator. 
Note also that the below quoted results are still valid in the more general context of
abstract {\it twin $BN$-pairs~} \cite[\S 1]{Abr97}. 

\begin{theorem}
\label{th - finiteness for parbolics}
Let $\Lambda$ be a Kac-Moody group over ${\bf F}_q$, and let $\Gamma$ be a facet fixator. 
Let $\Sigma$ be the Coxeter complex of the Weyl group $W$, i.e. the model for any
apartment in the building $X_\pm$ of $\Lambda$. 
Let $R$ be a chamber of $\Sigma$ and let $\Pi$ be the set of reflections in the codimension
one faces of $R$. 

\begin{enumerate}
\item[(i)] If any two reflections of $\Pi$ generate a finite group and if $q > 3$, then
$\Gamma$ is finitely generated; but it is not of cohomological type $FP_2$, in particular
not finitely presented, whenever some set of three reflections in $\Pi$ generates an infinite
subgroup of $W$. 
\item[(ii)] If any three reflections of $\Pi$ generate a finite group and if $q > 6$, then
$\Gamma$ is finitely presentable; but it is not of cohomological type $FP_3$ whenever some
set of four reflections in $\Pi$ generates an infinite subgroup of $W$. 
\end{enumerate}
\end{theorem}

When any two reflections of $\Pi$ generate a finite group, the Kac-Moody group is called 
{\it $2$-spherical}. 
Along with the Moufang property \cite[6.4]{RonLec} and twinnings \cite[Definition 3]{Abr97},
this notion plays a major role in the classification of buildings  with infinite Weyl
groups. 
In \cite{Abr97} some further results  are available; they deal with the higher
cohomological finiteness properties $FP_n$ and $F_n$ \cite[VIII.5]{BroCohomology}. 
According to Fact \ref{fact - finitely generated}, finite presentability is the
first finiteness property to be considered for the group $\Lambda$ itself. 
A result due to P. Abramenko and B. M\"uhlherr \cite{AbrMue} is available: 

\begin{theorem}
\label{th - finiteness for KM}
With the same notation as above, if any two reflections of $\Pi$ generate a finite group and
if $q > 3$, then $\Lambda$ is finitely presentable. 
\end{theorem}

For instance, finite presentability holds for Kac-Moody groups whose buildings have chambers
isomorphic to hyperbolic regular triangles of angle ${\pi \over 4}$ and $q>3$, but never holds
for those whose buildings are covered by chambers isomorphic to a regular right-angled
$r$-gon, $r \geq 5$ (we will see in \ref{ss - hyperbolic} that such Kac-Moody groups exist
in both cases). 

\pec 

\subsection{Continuous cohomology and Kazhdan's property (T)} 
\label{ss - T}
By results of J. Dymara and T. Januszkiewicz, being 2-spherical also implies useful
continuous cohomology vanishings  for automorphism groups of buildings. The result below is
a special case of \cite[Theorem E]{DJ02}. 

\begin{theorem}
\label{th - Kazhdan}
Let $\Lambda$ be a Kac-Moody group over $\mathbf{F}_q$, defined by a generalized 
Cartan matrix $A$ of size $n \times n$. 
Let $m<n$ be an integer such that all the principal submatrices of size $m \times m$ of $A$
are Cartan matrices (i.e. of finite type). 
Then for $1 \leq k \leq m-1$ and $q >\!\!> 1$, the continuous cohomology groups 
$H^k_{\rm ct}\bigl( {\rm Aut}(X_\pm),\rho \bigr)$ vanish for any unitary representation
$\rho$. 
\end{theorem}

The first cohomology case is extremely useful since vanishing of $H^1_{\rm ct}(G,\rho)$ for
any unitary representation $\rho$ is equivalent to property (T) \cite[Chap. 4]
{HarVal}. 
Therefore, when $\Lambda$ is 2-spherical, Theorem \ref{th - Kazhdan} implies property (T)
for the full automorphism groups
${\rm Aut}(X_\pm)$  with $q >\!\! >1$, hence for
their product, and finally for any lattice in this product, by S.P. Wang's Theorem 
\cite[III, Theorem 2.12]{Margulis}. 
The above result says in particular that many Kac-Moody groups have property (T), a
fundamental property for lattices of higher-rank simple algebraic groups
\cite[III]{Margulis}.

\subsection{Hyperbolic examples} 
\label{ss - hyperbolic}
The existence of buildings with prescribed shapes of apartements and links around
vertices is well-known in many cases \cite{RonFree}, \cite{BouMostow}. 
Some examples lead to interesting full automorphism groups and lattices \cite{BouPaj}, and
some other examples have surprisingly few automorphisms, though they are quite familiar
since they are Euclidean and tiled by regular triangles \cite{Bar}.
The result below \cite[Proposition 2.3]{RemGD} shows that some of these buildings are both
relevant to Kac-Moody theory and to (generalized) hyperbolic geometry \cite[III.H]{BriHae}. 
This enables to mix arguments of algebraic and geometric nature in the study of the
corresponding Kac-Moody groups. 

\pec 

\begin{prop}
\label{prop - hyperbolics exist} 
Let $P$ be a polyhedron in the hyperbolic space ${\Bbb H}^n$, with dihedral angles equal to
$0$ or to ${\pi \over m}$ with $m=2$, $3$, $4$ or $6$. 
For any prime power $q$, there is a Kac-Moody group whose building has constant
thickness equal to $q+1$ and where the apartments are all tilings of ${\Bbb H}^n$ by $P$;
this building is a complete geodesic ${\rm CAT}(-1)$-metric space. 
\end{prop}  

Recall that the {\it thickness~} at a codimension 1 cell $C$ of a building is the number of
maximal cells containing $C$. 
Note that according to G. Moussong \cite{Moussong}, any Coxeter
group acts discretely and cocompactly on a CAT(0)-space, which provides a good metric
realization of the associated Coxeter complex \cite[Chap. 2]{RonLec} (in which only
spherical facets are represented).  Moreover, the metric space under consideration is
CAT(-1) if and only if the Coxeter group is Gromov-hyperbolic (which is a weaker property in
general \cite[Chap. 3]{GhyHar}).  Using G. Moussong's non-positively curved complex as a
model for apartments, M. Davis
\cite{Davis} proved the existence of a CAT(0)-metric realization for any building, which is
CAT(-1) whenever the Weyl group of the building is Gromov-hyperbolic. 
In the latter case, the metric complex is usually not as nice as a hyperbolic tiling as in 
Proposition \ref{prop - hyperbolics exist}. 

\subsection{Easy non-linearities} 
\label{ss - easy NL}
On the one hand, in \ref{ss - lattice}-\ref{ss - T} we have quoted arguments
supporting the analogy between arithmetic groups over function fields and Kac-Moody groups
over finite fields.  
On the other hand, \ref{ss - hyperbolic} shows that the
geometries are certainly new since hyperbolic buildings can be obtained. 
Still, it is not proved at this stage that the groups are not well-known groups, only with
new actions: what about arguments showing that the groups are new? 
This question leads to rigidity problems, i.e. proving that groups acting naturally on some
geometries cannot act via a big quotient on another geometry. 
A way to attack the problem is to prove that the groups are not linear over any field. 
Easy non-linearities are summed up in the following \cite[Theorem 4.6]{RemNewton}: 

\begin{theorem}
\label{th - easy non-linearity} 
Let $\Lambda$ be a Kac-Moody group over $\mathbf{F}_q$  with infinite Weyl group, and let
$\Gamma$ be a facet fixator in $\Lambda$. 
Then $\Gamma$ always contains an infinite group of exponent $p$. 
Therefore $\Gamma$ cannot be linear over any field of characteristic $\neq p$. 
\end{theorem}

Exhibiting an infinite group of exponent $p$ is made possible thanks to arguments on
Kac-Moody root systems and pairs of parallel walls in infinite Coxeter complexes. 
Then, elementary Zariski closure and algebraic group arguments imply that an
infinite group of exponent $p$  cannot be linear over any field of characteristic different
from $p$ \cite[Lemma VIII.3.7]{Margulis}. 

\gec 

\section{Generalized Kac-Moody groups. Non-isomorphic groups with the same building}
\label{s - generalized KM}

In this section, we concentrate on buildings whose apartments are tilings of the
hyperbolic plane. 
Since we are interested in understanding new groups and geometries, it is natural to
consider the case of Fuchsian Weyl groups, corresponding to the  simplest exotic Kac-Moody
groups. 
The results are mainly of algebraic and combinatorial flavour; they emphasize the
role of conditions previously introduced to classify buildings. 
Abstract isomorphisms between Kac-Moody groups with the same hyperbolic building can be
naturally factorized (\ref{ss - automorphisms}). 
This is an analogy with the spherical building case, but a difference comes 
when exhibiting some non-isomorphic groups with the same building (\ref{ss -
non-isomorphic}). 
Moreover the local structure of right-angled Fuchsian buildings is
particularly simple: this  enables to construct groups which are very close to Kac-Moody
groups, but so to speak defined over several ground fields (\ref{ss - mixing}).  
Back to geometric group theory, this provides lattices of buildings of arbitrary large rank
satisfying a strong non-linearity property (\ref{ss - infinite kernel}). 

\subsection{Factorization of abstract automorphisms} 
\label{ss - automorphisms}
We say that a building is {\it Fuchsian~}if its Weyl group is the reflection
group of a tiling of the hyperbolic plane (in which case the latter tiling is a pleasant
model for apartments). 
According to Poincar\'e's Theorem \cite[4.H]{Maskit}, Fuchsian tilings are nice metric
realizations of Coxeter complexes (the most familiar ones after Euclidean tilings).  
Elementary facts from hyperbolic geometry enables to prove the following 
\cite[Theorem 3.1]{RemGD}: 

\begin{theorem}
\label{th - devissage} 
Let $G$ and $G'$ be two Kac-Moody defined over the same finite field ${\hbox {\bf F}_q}$ of 
cardinality $q \geq 4$. 
Assume that the associated buildings are all isomorphic $(\star)$ either to the same locally
finite tree, $(\star\star)$ or to the same Fuchsian building with regular chambers. 
Then, up to conjugacy in $G$, any abstract isomorphism from $G$ to $G'$ is the composition
of a permutation of the simple roots and possibly a global opposition of the sign of all
roots. 
\end{theorem}

This factorisation result is close to R. Steinberg's classical result on finite groups of Lie
type, saying that an abstract automorphism of such a group is the product of a Dynkin diagram
automorphism, a ground field automorphism and an inner automorphism \cite[Theorem
12.5.1]{Carter}. 
Other factorizations of isomorphisms have recently been obtained by P.-E. Caprace; the latter work is complementary since it deals with algebraically closed
fields \cite[Theorem 6]{Caprace}. 

\subsection{Several isomorphism classes} 
\label{ss - non-isomorphic}
Let $R$ be a right-angled $r$-gon of the hyperbolic plane ${\Bbb H}^2$.
For any integer $q \geq 2$ there is a unique building $I_{r,q+1}$ whose apartments are 
Poincar\'e tilings of ${\Bbb H}^2$ by $R$, and such that the link at each vertex is the
complete bipartite graph of parameters $(q+1,q+1)$ \cite[2.2.1]{BouMostow}. 
Here the {\it link~} of a vertex is a small ball around it, seen as a graph. 
We call such a building the {\it right-angled Fuchsian building~} of parameters $r$ and
$q+1$. 
These buildings are interesting because they locally look like products of trees, making 
them simple combinatorially, but globally their Weyl group is irreducible and
Fuchsian, leading to remarkable rigidity properties \cite{BouMostow}, \cite{BouPaj}. 
By uniqueness and Proposition \ref{prop - hyperbolics exist}, for each $r \geq 5$
the building $I_{r,q+1}$ comes from a Kac-Moody group whenever $q$ is a prime power. 
Using Theorem \ref{th - devissage}, it can be proved that there are abstractly
non-isomorphic Kac-Moody groups with the same building \cite[\S 4, Proposition]{RemGD}: 

\begin{cor}
\label{cor - non-isomorphic} 
Let $q \geq 4$ be a prime power and let $r \geq 5$ be an integer. 
Then there are several abstract isomorphism classes of Kac-Moody groups whose associated
buildings are the same $I_{r,q+1}$. 
\end{cor}

This result is in contrast with the spherical case, where according to J. Tits'
classification, a spherical building of rank $\geq 3$ uniquely determines a field and 
an algebraic group over this field \cite{TitsSpherical}. 
Here the rank $r$ of the building is an arbitrary integer $\geq 5$. 
This is an argument explaining why in the classification of Moufang
twinnings \cite{MuhRon}, \cite{Muhlherr}, the buildings must be assumed 2-spherical: 
Kac-Moody groups with $I_{r,q+1}$ buildings obviously do not satisfy this property. 

\subsection{Mixing ground fields} 
\label{ss - mixing}
There is an even stronger argument showing that being 2-spherical is a necessary condition
for a Moufang twin building to be part of a reasonable classification (where the group side
would be given by Kac-Moody groups and their twisted versions). 
Indeed, some generalizations of Kac-Moody groups with several ground fields can be
constructed, provided the buildings have apartments tiled by regular hyperbolic
right-angled $r$-gons. 
Here \og generalization\fg means that the groups satisfy the same
combinatorial axioms as Kac-Moody groups (namely, those of twin root data \cite{Tit92}).  
We have \cite[Theorem 3.E]{RemRon}: 

\begin{theorem}
\label{th - exotic Fuchsian}
A right-angled Fuchsian building belongs to a Moufang twinning whenever its thicknesses
at panels are cardinalities of projective lines.
\end{theorem}

In the case of trees, a more theoretical construction is sketched by J. Tits in his Notes de
Cours au Coll\`ege de France \cite[\S 9]{TitsCDF}. 
Still, the rank of the infinite dihedral Weyl group there is always equal to 2, whereas in
the above result it can be any integer $\geq 5$. 

\subsection{Stronger non-linearities} 
\label{ss - infinite kernel}
Besides its combinatorial interest, the above construction can be seen as a way to produce
lattices of hyperbolic buildings with remarkable non-linearity properties. 
Recall that according to Theorem \ref{th - easy non-linearity} the characteristic $p$ of the
ground field of a Kac-Moody group prevents this group from being linear over any field of
characteristic $\neq p$. 
Therefore, mixing fields of different characteristics in Theorem \ref{ss - mixing} should
enable to produce groups which are not linear over any field. 
This is indeed the case \cite[Theorem 4.A]{RemRon}: 

\begin{theorem}
\label{th - infinite kernel}
Let $r\geq 5$ be an integer and $\{ {\bf K}_i \}_{i \in {\bf Z}/r}$ be a family of 
fields, among which are two fields with different positive characteristics. 
Let $\Lambda$ be a group defined as in Theorem \ref{th - exotic Fuchsian} from these fields
and let $\Gamma$ be a chamber fixator. 
Then, any group homomorphism $\rho: \Gamma \to \prod_{\alpha \in A}  {\bf G}_\alpha({\bf
F}_\alpha)$ has infinite kernel, whenever the index set $A$ is finite and ${\bf  G}_\alpha$
is a linear algebraic group defined over the field ${\bf F}_\alpha$ for each $\alpha \! \in
\! A$.
\end{theorem}

Note that the result is stronger than a plain non-linearity since mixing ground fields is
also allowed at the (linear) right hand-side of the representation. 
Moreover the kernel is not only non-trivial, but always infinite. 
From the purely Kac-Moody viewpoint, this result says that the truly difficult case of
non-linearity is when the characteristic of the target algebraic group is the characteristic
of the ground field of the Kac-Moofy group. 
This case is reviewed in Sect. \ref{s - hard NL}. 

\gec 

\section{Generalized algebraic groups over local fields}
\label{s - topological}

A Kac-Moody group $\Lambda$ acts discretely on the product of its buildings, but its action
on a single factor is no longer discrete. 
Therefore it makes sense to take the closure $\overline\Lambda < {\rm Aut}(X_\pm)$ of the
image of a Kac-Moody group in such a non-discrete action. 
The result is called a {\it topological Kac-Moody group~} (the kernel of the
$\Lambda$-action on $X_\pm$ is the finite center of $\Lambda$). 
In the classical case $\Lambda={\rm SL}_n({\bf F}_q[t,t^{-1}])$, $X_\pm$ is the
Bruhat-Tits building of ${\rm SL}_n \bigl( {\bf F}_q(\!(t)\!) \bigr)$ or 
${\rm SL}_n \bigl( {\bf F}_q(\!(t^{-1})\!) \bigr)$, respectively. 
If $\mu_n({\bf F}_q)$ denotes the $n$-th roots of unity in ${\bf F}_q$, the image
$\Lambda/Z(\Lambda)$ of ${\rm SL}_n({\bf F}_q[t,t^{-1}])$ under the action on $X_\pm$ is
${\rm  SL}_n({\bf F}_q[t,t^{-1}])/\mu_n({\bf F}_q)$ and the
completions $\overline\Lambda$ are ${\rm PSL}_n \bigl( {\bf F}_q(\!(t)\!) \bigr)$ and 
${\rm PSL}_n \bigl( {\bf F}_q(\!(t^{-1})\!) \bigr)$, respectively. 
In fact, there are many arguments to
compare a topological Kac-Moody group with a semisimple group over a local field:
existence of a combinatorial structure refining Tits systems (\ref{ss - RTS}),
group-theoretic characterization of chamber-fixators in terms of pro-$p$ subgroups 
(\ref{ss - Iwahori}), topological simplicity (\ref{ss - top simple}). 
Still, as for discrete groups we have phenomena suggesting that some so-obtained
totally disconnected groups are new (\ref{ss - F boundaries}). 
We finally quote a result exhibiting the coexistence of non-linear non-uniform lattices and
uniform lattices embedding convex-cocompactly in hyperbolic spaces (\ref{ss - coexistence}). 

\subsection{Refined Tits systems} 
\label{ss - RTS}
The structure of a {\it refined Tits system~} is due to V. Kac and D. Peterson
\cite{KacPet}; it is a generalization of a split $BN$-pair, a notion from the theory of finite
groups of Lie type.  The difference with a plain Tits system is the formalization in the
axioms of the existence of an abstract unipotent subgroup in the Borel subgroup. 
We have \cite[Theorem 1.C]{RemRon}: 

\begin{theorem}
\label{th - refined TS}
Let $\Lambda$ be a Kac-Moody group over the finite field ${\bf F}_q$ of characteristic $p$. 
The associated topological Kac-Moody group $\overline\Lambda$ admits a refined Tits system,
and the Tits system gives rise to a building in which any facet-fixator is, up to finite index,
a pro-$p$ group.
\end{theorem}

It is well-known that a group acting strongly transitively on a building 
admits a natural Tits system \cite[Theorem 5.2]{RonLec}, so the point in the first
assertion lies in the difference between a Tits system and a refined Tits
system. 
Moreover standard properties of double cosets in Bruhat decompositions imply
that topological Kac-Moody group are compactly generated \cite[Corollary 1.B.1]{RemLin}.
It is explained in \cite[1.B.1]{RemLin} why refined Tits systems for $\overline\Lambda$
imply the existence of a lot of torsion in facet fixators -- called {\it parahoric
subgroups}. 
This explains why the analogy with semisimple groups over local fields is relevant
only when the local field has the same characteristic $p$ as the finite ground field of
$\Lambda$ (parahoric subgroups of semisimple groups over finite extensions of ${\bf Q}_l$
are virtually torsion free). 
In analogy with the classical case of Bruhat-Tits theory \cite{BruhatTits1},
\cite{BruhatTits2}, we call {\it Iwahori subgroup~} a chamber fixator in 
$\overline\Lambda$.  Ê

\subsection{Iwahori subgroups} 
\label{ss - Iwahori}
Back to the case $\Lambda={\rm SL}_n({\bf F}_q[t,t^{-1}])$, let $v$ be a vertex in the
Bruhat-Tits building of ${\rm SL}_n \bigl(  {\bf F}_q(\!(t)\!) \bigr)$.
Then its fixator is isomorphic to ${\rm SL}_n({\bf F}_q[[t]])$, and for some chamber
containing $v$ the corresponding Iwahori subgroup is the group of matrices in 
${\rm SL}_n({\bf F}_q[[t]])$ reducing to upper triangular matrices modulo $t$. 
Moreover the first congruence subgroup of ${\rm SL}_n({\bf F}_q[[t]])$, i.e. the matrices
reducing to the identity modulo $t$, is a maximal pro-$p$ subgroup whose normalizer is the 
above Iwahori subgroup. 
In the general Kac-Moody setting, the result below \cite[Proposition 1.B.2]{RemLin} is a
generalization of \cite[Theorem 3.10]{PR94}.

\begin{prop}
\label{prop - Iwahori} 
The Iwahori subgroups are group-theoretically characterized as the normalizers of
the  pro-$p$ Sylow subgroups.
\end{prop}

In general, the analogue of the first congruence subgroup of a vertex fixator must be
defined as the (pointwise) fixator of the star around the vertex under consideration. 

\subsection{Topological simplicity} 
\label{ss - top simple}
In view of the simplicity of adjoint algebraic simple groups over large enough fields, the
theorem below \cite[Theorem 2.A.1]{RemLin} is not surprising: 

\begin{theorem}
\label{th - top simple}
For $q \geq 4$, a topological Kac-Moody group over a finite field is the direct product of
finitely many topologically simple groups, with one factor for each connected  component of
its Dynkin diagram.
\end{theorem}

This result is extremely useful when extending abstract representations of finitely 
generated Kac-Moody to continuous representations of topological Kac-Moody groups 
(Theorem \ref{th - embedding}).  
The arguments of the proof are basically: a normal subgroup in an irreducible Tits system
either is chamber-transitive or acts trivially on the associated building
\cite[IV.2]{Bou81},  Iwahori subgroups are virtually pro-$p$, and a generating set for a
Kac-Moody group can be chosen in a finite collection of finite subgroups of Lie type (which
are perfect whenever $q \geq 4$). 
It is reasonable to think that topological Kac-Moody groups are in fact abstractly simple
(Question \ref{quest - abstractly simple}). 

\subsection{Non-homogeneous Furstenberg boundaries} 
\label{ss - F boundaries}
When the building $X_\pm$ has hyperbolic apartments, the existence of many hyperbolic
translations leads to an interesting connection with topological dynamics. 
Let $Y$ be a compact metrizable space and let us denote by $M^1(Y)$ the space
of probability measures on $Y$. 
Recall that $M^1(Y)$ is compact and metrizable for the weak-$*$ topology.
If $Y$ admits a continuous action by a topological group $G$, we say that $Y$ is a {\it
Furstenberg boundary~} for $G$ if it is $G$-{\it minimal~} and $G$-{\it strongly proximal~} 
\cite[VI.1.5]{Margulis}.  
The first condition says that any $G$-orbit in $Y$ is dense and the second one says that any 
$G$-orbit closure in $M^1(Y)$ contains a Dirac mass. 
The arguments of \cite[Lemma 4.B.1]{RemLin} give: 

\pec

\begin{lemma}
\label{lemma - Furstenberg}
Let $\Lambda$ be a Kac-Moody group whose buildings have apartments isomorphic to a
hyperbolic tiling. 
Then the asymptotic boundary $\partial_\infty X_\pm$ is a Furstenberg boundary for any closed
automorphism group of $X_\pm$ containing $\overline\Lambda$.
\end{lemma}

This existence of non-homogeneous boundaries is new with respect to semisimple Lie groups 
(archimedean or not), since in the latter case any Furstenberg boundary is equivariantly
isomorphic to a flag variety of the group \cite[9.37]{GJT}.
	
\subsection{Coexistence of two kinds of lattices} 
\label{ss - coexistence} 
We close this section by stating a result which says that some topological Kac-Moody
groups contain lattices of surprisingly different nature \cite[Proposition 4.B]{RemRon}. 

\begin{prop}
\label{prop - two lattices}
There exist topological Kac-Moody groups $\overline\Lambda$ over ${\bf F}_q$ which contain
both non-uniform lattices which cannot be linear over any field of 
characteristic prime to $q$,
and uniform lattices which have convex-cocompact embeddings into real 
hyperbolic spaces.
The limit sets of the latter embeddings often have Hausdorff dimension $>2$.
\end{prop}

It follows from \cite[Corollary 0.5]{BurMoz96} that if  a Kac-Moody group $\Lambda$ is
$S$-arithmetic and such that $\overline\Lambda$ is a higher-rank simple Lie group, then any
lattice of $\overline\Lambda$ fixes a point in each of its actions on proper
CAT(-1)-spaces.  Therefore the above phenomenon is excluded in the classical algebraic
situation, unless the building of $\overline\Lambda$ is a tree, meaning that the Weyl group
$\Lambda$ is infinite dihedral (of rank 2 as a Coxeter group). 
The rank $r$ of a Fuchsian building $I_{r,q+1}$ may be chosen arbitrarily large $\geq 5$. 

\gec 

\section{Non-linearity in equal characteristics}
\label{s - hard NL} 

In view of the easy non-linearity results (Theorem \ref{th - easy non-linearity}), the
remaining linearity to disprove is for Kac-Moody groups over finite fields of characteristic
$p$ into linear groups of characteristic $p$. 
But some Kac-Moody $S$-arithmetic groups such as ${\rm SL}_n \bigl( {\bf F}_q[t,t^{-1}]
\bigr)$  are linear in equal characteristic. 
Therefore the question is rather, for any finite field of characteristic $p$, to find
examples of Kac-Moody groups that cannot be linear over any field of characteristic $p$
either.  
The main idea is to use some superrigidity property (\ref{ss - commensurator
superrigidity}) to show that the existence of a faithful abstract homomorphism from a
finitely generated Kac-Moody group into an algebraic group implies the existence of
an embedding of a topological Kac-Moody group into a non-Archimedean simple group 
(\ref{ss - embedding}).  
The existence of such an embedding is expected to be simpler to disprove because it comes
with an embedding of the vertices of the (possibly exotic) Kac-Moody building to a Euclidean
one, which enables to take advantage of the incompatibility between hyperbolic and Euclidean
geometries.  For some Kac-Moody groups with Fuchsian buildings this is indeed the case 
(\ref{ss - NL Fuchsian}). 
The same circle of ideas allows to show a complementary result: a non-faithful representation
from a finitely generated Kac-Moody group most of the time has a virtually solvable image 
(\ref{ss - v solvable}). 
Combining the latter results leads to groups all of whose linear images are
finite (\ref{ss - finite linear images}). 

\pec 

\subsection{Commensurator superrigidity} 
\label{ss - commensurator superrigidity}
Let us first recall that the {\it commensurator~} of a group inclusion $\Delta<G$ is the
group:

\pec 

\centerline{
${\rm Comm}_{G}(\Delta):= \{ g \! \in \! G \mid \Delta \cap 
g\Delta g^{-1} \hbox{\rm ~has finite
index in both } \Delta \hbox{\rm ~and } g \Delta g^{-1}\}$.} 

\pec 

The next theorem is basically due to G.A. Margulis \cite[VII.5.4]{Margulis}, and the
formulation below can be found in \cite[Theorem 1]{Bonvin}. 
G.A. Margulis proved it when $G$ is a semisimple group over a local field. 
A certainly non-exhaustive list of later contributions is the following: 
in \cite{AcaBur}, the semisimple group $G$ is replaced by a group containing an amenable 
subgroup $P<G$ similar to a minimal parabolic subgroup; 
in \cite{Burger} the existence of $P$ is replaced by the existence of a substitute for a 
Furstenberg boundary; 
in \cite{BurMon} such boundaries are constructed for compactly generated groups $G$; 
and the double ergodicity Theorem in \cite{Kaim} shows that suitable boundaries are
available for any locally compact second countable group, via Poisson boundary theory. 
Different approaches lead to similar results: \cite{Marg94} by means of equivariant
generalized harmonic mappings and \cite{Shalom} by means of representation theory. 

\begin{theorem} 
\label{th - superrigidity}
Let $G$ be a locally compact second countable topological group, $\Gamma < G$ be a lattice
and $\Lambda$ be a subgroup of $G$ with $\Gamma < \Lambda < {\rm Comm}_G (\Gamma)$.
Let $k$ be a local field and $H$ be a connected almost $k$-simple algebraic group.
Assume $\pi: \Lambda \to H_k$ is a homomorphism such that $\pi(\Lambda)$ is dense in the
Zariski topology on $H$ and $\pi(\Gamma)$ is unbounded in the Hausdorff topology on $H_k$.
Then $\pi$ extends to a continuous homomorphism $\overline{\Lambda} \to H_k / Z(H_k)$,
where $Z(H_k)$ is the center of $H_k$.
\end{theorem}

A more difficult superrigidity consists in extending representations of lattices
(instead of their commensurators) to continuous representations of the ambient
topological group (to algebraic representations when so is the ambient group).
This is a more difficult result which requires stronger assumptions (e.g. higher rank
for algebraic ambient groups), and again the main ideas are due to G.A. Margulis.  
The first results in positive characteristic were proved by T.N. Venkataramana
\cite{Venka}. 

\pec 

\subsection{Embedding theorem} 
\label{ss - embedding}
For the following theorem, we were inspired by a paper due to A. Lubotzky, Sh. Mozes and
R.J. Zimmer \cite{LMZ94}, where superrigidity is used to disprove
linearities for commensurators of tree lattices. 
Trees are one-dimensional buildings, and Kac-Moody theory naturally leads to lattices for
buildings of arbitrary dimension. 
We have \cite[1.B Lemma 2]{RemRon}: 

\begin{lemma} 
\label{lemma - commensurator}
Let $\Lambda$ be any Kac-Moody group over ${\bf F}_q$ and let $\Gamma$ be any negative facet
fixator in $\Lambda$. 
Then: $\Lambda < {\rm Comm}_{\overline\Lambda}(\Gamma)$. 
\end{lemma}

Therefore it makes sense to use commensurator superrigidity because when $q >\!\!> 1$,  
$\Gamma$ is a lattice of $\overline\Lambda$ by Theorem  \ref{th - lattice}. 
This enables to obtain \cite[\S 3, Theorem]{RemLin}: 

\begin{theorem}
\label{th - embedding}
Let $\Lambda$ be a Kac-Moody group over the finite 
field ${\bf F}_q$ of
characteristic $p$ with $q > 4$ elements, with infinite Weyl group 
$W$ and buildings $X_+$ and
$X_-$.
Let $\overline\Lambda$ be the corresponding Kac-Moody topological group.
We make the following assumptions:

\pec

\hskip 4mm {\rm (TS)~} the group $\overline\Lambda$ is topologically simple;

\hskip 4mm {\rm (NA)~} the group $\overline\Lambda$ is not amenable;

\hskip 4mm {\rm (LT)~} the group $\Lambda$ is a lattice of $X_+ \times X_-$ for its diagonal
action.

\pec

Then, if $\Lambda$ is linear over a field of characteristic $p$, there exists:

-- a local field $k$ of characteristic $p$ and a connected adjoint $k$-simple
group ${\bf G}$,

-- a topological embedding $\mu: \overline\Lambda \to {\bf G}(k)$ 
with  Hausdorff unbounded
and Zariski dense image,

-- and a $\mu$-equivariant embedding $\iota: V_{X_+} \to V_\Delta$
from the set of vertices of the Kac-Moody building $X_+$ of $\Lambda$ 
to the set of vertices of
the Bruhat-Tits building $\Delta$ of ${\bf G}(k)$.
\end{theorem}

Conditions (TS) and (LT) are satisfied whenever the Weyl group $W$ of
$\Lambda$ is irreducible and $q >\!\!> 1$. 
Condition (NA) is discussed in \ref{s - conjectures}, where it is conjectured that
it is empty (Problem \ref{pb - NA}). 
Note that the conclusion of Theorem \ref{th - superrigidity} provides a
continuous extension but it is not clear that this topological homomorphism is a closed map.
Indeed, there is still work to do \cite[Lemma 3.C]{RemLin}
and, besides the topological simplicity of
$\overline\Lambda$, the arguments used for that are of combinatorial nature: basically
the Bruhat decomposition of $\overline\Lambda$ with respect to an Iwahori subgroup. 

\subsection{Non-linear Fuchsian Kac-Moody groups} 
\label{ss - NL Fuchsian} 
In the case of some Kac-Moody groups with Fuchsian buildings, it can be proved that the
associated topological groups cannot be closed subgroups of any simple non-Archimedean Lie
group, thus implying the non-linearity of the involved finitely generated groups
\cite[Theorem 4.C.1]{RemLin}. 

\begin{theorem}
\label{th - Fuchsian NL}
Let $\Lambda$ be a countable Kac-Moody group over ${\bf F}_q$ with right-angled
Fuchsian associated buildings. 
Assume that any prenilpotent pair of roots not contained in a spherical root system leads 
to a trivial commutation of the corresponding root groups. 
Then for $q >\!\!> 1$, the group $\Lambda$ is not linear over any field. 
\end{theorem}

The condition on commutation of root groups is technical: prenilpotency of pairs
of roots is relevant to abstract root systems of Coxeter groups \cite{Tit87},
\cite[1.4.1]{RemAst}. 
It is not very restrictive and actually a weaker assumption may be
required -- see the remark after \cite[Lemma 4.A.2]{RemLin}. 
Seeing roots as half-apartments, a pair of roots is prenilpotent if and only if the
walls of the roots intersect or if a root contains the other. 

\pec 

The main idea at this stage is to use a dynamical characterization of parabolic subgroups
proved by G. Prasad in a paper on strong approximation \cite[Lemma 2.4]{Prasad} in the
algebraic case,  e.g. the case of the target ${\bf G}(k)$ of the continuous extension 
$\mu : \overline\Lambda \to {\bf G}(k)$. 
In the case of hyperbolic buildings, one must first say that a parabolic subgroup is by
definition a boundary point fixator (in analogy with the symmetric space or Bruhat-Tits
building case), and then show that the dynamical characterization holds too \cite[Lemma
4.B.2]{RemLin}. 
We now have a dynamical roundabout to the non-existence of an algebraic structure on
$\overline\Lambda$, and this enables to show that under the continuous homomorphism $\mu$, 
parabolic subgroups go to parabolic subgroups. 
The contradiction comes when one notes that, thanks to the strong dynamics of a Fuchsian
group on the asymptotic boundary of the hyperbolic plane, the dynamical analogues of
unipotent radicals on the left hand-side are not normalized by the parabolic subgroups
containing them, whereas it is the case by definition on the right 
hand-side of $\mu$.

\pec

\subsection{Virtual solvability of non-faithful linear images}
\label{ss - v solvable}
In \ref{ss - embedding}, the starting point is a faithful abstract representation from a
finitely generated Kac-Moody group. 
Starting from non-faithful representations is interesting too, since it can be shown that in
this case the images are solvable up to finite index. 
In other words, when the kernel is non-trivial, it is big \cite[Theorem 11]{RemImage}: 

\begin{theorem}
\label{th - virtually solvabe}
Let $\Lambda$ be a Kac-Moody group over the finite field ${\bf F}_q$ of characteristic $p$, 
with connected Dynkin diagram and $q >\!\!> 1$. 
Let $\rho: \Lambda \to {\rm GL}_n({\bf F})$ be a linear representation. 
If $\rho$ is not faithful, the group $\rho(\Lambda)$ is virtually solvable. 
In particular if $\Lambda$ is Kazhdan, $\rho(\Lambda)$ is finite. 
\end{theorem}

Recall that by Theorem \ref{th - Kazhdan}, many Kac-Moody groups over ${\bf F}_q$ are Kazhdan
whenever $q >\!\!> 1$. 
The above theorem is proved thanks to the same kind of arguments as for Theorem \ref{th -
embedding}. 
In order to be in position to apply Theorem \ref{th - superrigidity}, we have to take the
Zariski closure of the image $\rho(\Lambda)$ and to mod out by the radical 
$R\bigl(\overline{\rho(\Lambda)}^Z\bigr)$ of the latter algebraic group. 
It then suffices to show that the image of $\Lambda$ in 
$\overline{\rho(\Lambda)}^Z/R\bigl(\overline{\rho(\Lambda)}^Z\bigr)$ is finite,
which can be done thanks to Burnside's theorem \cite[4.5, Exercise 8]{Jacobson}. 
A way to fulfill the unboundedness assumption on the image of $\Lambda$ (in the Hausdorff
topology) is to use tricks from Tits' paper on the existence of free groups in linear groups
\cite{Tit72}.

\pec 

\subsection{Groups all of whose linear images are finite} 
\label{ss - finite linear images} 
From the previous result, it can be shown that some Kac-Moody groups with hyperbolic
buildings don't have any infinite linear image \cite[Theorem 16]{RemImage}. 

\begin{theorem}
\label{th - finite linear images}
There is an integer $N$ such that for any Kac-Moody group $\Lambda$ over ${\bf F}_q$ whose
buildings have apartments isomorphic to the tesselation of the hyperbolic plane by regular
triangles of angle  ${\pi \over 4}$, if $q \geq N$ then any linear image of $\Lambda$ is
finite, whatever the target field. 
\end{theorem}

The proof roughly goes as follows. 
By Theorem \ref{th - virtually solvabe}, it is enough to show that some Kac-Moody groups
enjoying Kazhdan's property (T) are not linear over any field. 
According to Theorem \ref{th - Kazhdan}, the groups as in the theorem are
Kazhdan, and by twin root datum arguments it can be proved that such groups contain
Kac-Moody-like subgroups to which the proof of Theorem \ref{th - Fuchsian NL} applies. 

\gec

\section{Conjectures and questions}
\label{s - conjectures} 

In this final section, we propose a few questions about finitely generated and 
totally disconnected Kac-Moody groups. 
Here, the understatement in {\it totally disconnected~}Êis that the group under
consideration is uncountable and not endowed with the discrete topology. 

\subsection{Finitely generated groups} 
\label{ss - q discrete} 
For discrete groups, we think that the class of non-linear Kac-Moody groups
is much wider than the one for which the property has been proved so far. 

\begin{conj}
\label{conj - hyperbolic implies non-linear}
If the Weyl group $W$ is Gromov-hyperbolic and if $q >\!\!> 1$, then the 
Kac-Moody group $\Lambda$ is not linear. 
\end{conj}
 
Recall that according to G. Moussong, in the class of Coxeter groups hyperbolicity is
equivalent to acting discretely and cocompactly on a ${\rm CAT}(-1)$-space \cite{Moussong}. 
The non-linearity proof of Sect. \ref{s - hard NL} deals with groups $\Lambda$ whose Weyl
groups have pleasant Coxeter complexes since the latter complexes are Fuchsian tilings. 
In general, Moussong's complex is defined by isometric gluings of cells and is less easy to
understand, so one may have to use arguments of combinatorial nature on Kac-Moody
root systems at some points. 
Note also that one should get rid of the technical condition on prenilpotent pairs of roots
in the statement of Theorem \ref{th - Fuchsian NL}. 

\pec 

Another natural question is: 

\begin{quest}
\label{quest - on q}
Can the assumption \og $q >\!\!> 1$\fg be removed in non-linearity results? 
\end{quest}

In other words: is there a generalized Cartan matrix $A$ and a prime number $p$ such that the
corresponding group is linear over ${\bf Z}/p$ but no longer linear over 
${\bf F}_{p^r}$ for large enough $r$? 
Note that for $\Lambda$ to be a lattice of the product of its buildings -- which is a
crucial argument in our proof, the value of $q$ is important (Theorem \ref{th - lattice}). 
So if the answer to the above problem is yes, the non-linearity proof should require new
ideas. 
Moreover a counter-example due to P. Abramenko shows that for buildings whose chambers are
regular hyperbolic triangles of  angles ${\pi \over 4}$, facet fixators in some groups over
${\bf F}_2$ or ${\bf F}_3$ are not finitely generated \cite[Counter-example 1, Remark
2]{AbrBie}, hence cannot have property (T) \cite[Theorem III.2.7]{Margulis}. 
But according to Theorem \ref{th - Kazhdan}, for $q >\!\!> 1$ the groups do enjoy Kazhdan's
property (T). 

\pec 

The most general question on non-linearity for Kac-Moody groups is: 

\begin{pb}
\label{pb - GCM}
Find necessary and sufficient conditions for the non-linearity of a finitely generated
Kac-Moody group, only in terms of the generalized Cartan matrix defining the group. 
\end{pb}

Of course, if the generalized Cartan matrix $A$ is of affine type, meaning that the Weyl
group is a Euclidean reflection group, then the corresponding Kac-Moody groups are
$S$-arithmetic groups, hence are linear. 
A restatement is: are there other linear examples than the affine ones? 

\pec 

At last, there is another problem which comes from the analogy with the situation of
lattices, namely the problem of arithmeticity. 
Of course when algebraic structures are available, the definition is well-known 
\cite[Definitions 6.1.1 and 10.1.11]{Zimmer}. 
In the general context of an inclusion $\Delta < G$ of a discrete subgroup $\Delta$ in a
locally compact group $G$, it has now become classical to say that $\Delta$ is {\it
arithmetic~}Êin $G$ if by definition ${\rm Comm}_G(\Delta)$ is dense in $G$ 
(\ref{ss - commensurator superrigidity}). 
This is reasonable since by Margulis' commensurator criterion, a lattice in a non-compact
simple Lie group is arithmetic in the classical sense if and only if so it is in the
above sense \cite[Theorem 6.2.5]{Zimmer}. 
This definition makes sense in the case where the topological group $G$ is the isometry group
of some metric space, e.g. of some building with enough autmorphisms. 
Since Kac-Moody buildings are close with many respects to Bruhat-Tits buildings, it is
natural to ask: 

\begin{quest}
\label{quest - arithmetic}
Let $X_+$ be the positive building of some Kac-Moody group $\Lambda$ over ${\bf F}_q$. 
Is a negative facet fixator in $\Lambda$ arithmetic in ${\rm Aut}(X_+)$? 
More generally, which lattices, not necessarily relevant to Kac-Moody groups, are arithmetic
in ${\rm Aut}(X_+)$? 
\end{quest}

Note that by definition of a topological Kac-Moody group (Sect. \ref{s - topological},
introduction), any negative facet fixator $\Gamma$ is arithmetic in $\overline\Lambda$ since
by Lemma \ref{lemma - commensurator}, ${\rm Comm}_{\overline\Lambda}(\Gamma)$ contains
$\Lambda$. 
Here is a list of known arithmeticity results: 

\pec 

\begin{enumerate}

\item[(i)] Cocompact lattices in a locally finite tree $T$ are all arithmetic in ${\rm
Aut}(T)$ \cite{Liu}. 

\item[(ii)] The Nagao lattice ${\rm SL}_2 \bigl( {\bf F}_q [t^{-1}] \bigr)$
is a non-uniform arithmetic lattice in the full automorphism group of the Bruhat-Tits tree of
${\rm SL}_2 \bigl( {\bf F}_q (\!( t )\!) \bigr)$ \cite{MozNagao}. 
This was later extended to some Moufang twin trees, which enables to deal with non-linear
non-uniform tree lattices \cite{AbrRem}. 

\item[(iii)] Uniform lattices in some hyperbolic buildings, not necessarily from Kac-Moody
theory \cite{HagComm}. 
\end{enumerate}

See the end of \cite{AbrRem} for more precise questions. 

\subsection{Totally disconnected groups} 
\label{ss - q td} 
Theorem \ref{th - embedding} shows that non-linearity of some finitely generated groups
roughly amounts to non-linearity of much bigger topological groups. 
The study of topological Kac-Moody groups in Sect. \ref{s - topological} made appear an
interesting class of pro-$p$ subgroups. 
A way to disprove some linearities would be to answer the following: 

\begin{quest}
\label{quest - non-linear pro-p}
Which topological Kac-Moody groups contain non-linear pro-$p$ subgroups? 
\end{quest}

Forgetting discrete groups, the question is interesting in its own: for instance, the
linearity of free pro-$p$ subgroups is a question so far admitting only partial answers 
(E. Zelmanov).  

\pec 

In the classical case, topological Kac-Moody groups correspond to groups of rational points
of adjoint semisimple groups, so it makes sense to ask the following: 

\begin{quest}
\label{quest - abstractly simple}
Are topological Kac-Moody groups direct products of abstractly simple groups? 
\end{quest}

Theorem \ref{th - top simple} only says that the latter groups are direct products
of topologically simple groups. 
F. Haglund and F. Paulin, elaborating on J. Tits' proof for trees \cite{TitsTree},
showed the abstract simplicity of full automorphism groups of many hyperbolic buildings 
\cite[Theorems 1.1 and 1.2]{HagPau}. 

\pec 

At last, in two important cases topological Kac-Moody groups $\overline\Lambda$ are not
amenable. 
(Recall that a simple Lie group, when non-compact, is not amenable.) 

\begin{pb}
\label{pb - NA}
Show that topological Kac-Moody groups are never amenable. 
\end{pb}

Solving the problem, even for $q>\!\!>1$, would enable to remove assumption (NA) in Theorem
\ref{th - embedding}. 
It is solved when the Weyl group $W$ of $\Lambda$ is Gromov-hyperbolic, and when property
(T) holds for $\overline\Lambda$. 
The latter case is clear because amenability and property (T) imply compactness
\cite[Corollary 7.1.9]{Zimmer}, which then implies to have a fixed point in the building
$X_+$ by non-positive curvature \cite[VI.4]{BroBuildings}: contradiction with the
chamber-transitivity of the $\Lambda$-action on $X_+$. 
The former case follows from a Furstenberg lemma for CAT(-1)-spaces \cite[Lemma
2.3]{BurMoz96} -- see \cite[Sect. 3, Introduction]{RemLin} for details. 

\bibliographystyle{amsalpha}
\bibliography{DiscreteSurvey}

\providecommand{\bysame}{\leavevmode\hbox to3em{\hrulefill}\thinspace}
\begin{thebibliography}{Rem03b}

\bibitem[AB94]{AcaBur}
N.~A'Campo and M.~Burger, \emph{R\'eseaux arithm\'etiques et commensurateur
  d'apr\`es {G}.{A}. {M}argulis}, Inventiones Mathematic\ae\, \textbf{116}
  (1994), 1--25.

\bibitem[Abr97]{Abr97}
P.~Abramenko, \emph{Twin buildings and applications to {$S$}-arithmetic
  groups}, Lecture Notes in Mathematics, vol. 1641, Springer, 1997.

\bibitem[Abr03]{AbrBie}
P.~Abramenko, \emph{Finiteness properties for groups acting on twin buildings},
  Groups: Topological, Combinatorial and Arithmetic Aspects (Bielefeld, 1999)
  (T.W. M\"uller, ed.), LMS Lecture Notes Series, London Mathematical Society,
  Cambridge University Press, 2003, pp.~17--21.

\bibitem[AM97]{AbrMue}
P.~Abramenko and B.~M\"uhlherr, \emph{Pr\'esentations de certaines {BN}-paires
  jumel\'ees comme sommes amalgam\'ees}, Comptes-Rendus de l'Acad\'emie des
  sciences de Paris \textbf{325} (1997), 701--706.

\bibitem[AR03]{AbrRem}
P.~Abramenko and B.~Remy, \emph{Commensurators of some nonuniform tree lattices
  and {M}oufang twin trees}, this volume, 2003.

\bibitem[Bar00]{Bar}
S.~Barr\'e, \emph{Immeubles de {T}its triangulaires exotiques}, Ann. Fac. Sci.
  Toulouse Math. \textbf{9} (2000), no.~4, 575--603.

\bibitem[Beh69]{BehrRed}
H.~Behr, \emph{Endliche {E}rzeugbarkeit arithmetischer {G}ruppen \"uber
  {F}unktionenk\"orpern}, Inventiones Mathematic\ae\, \textbf{7} (1969), 1--32.

\bibitem[Beh03]{Behr}
H.~Behr, \emph{Higher finiteness properties for {S}-arithmetic groups in the
  function field case {I}}, Groups: Topological, Combinatorial and Arithmetic
  Aspects (Bielefeld, 1999) (T.W. M\"uller, ed.), LMS Lecture Notes Series,
  London Mathematical Society, Cambridge University Press, 2003, pp.~761--769.

\bibitem[BH99]{BriHae}
M.~Bridson and A.~Haefliger, \emph{Metric spaces of non-positive curvature},
  Grundleheren der mathematischen Wissenschaften, vol. 319, Springer, 1999.

\bibitem[BM96]{BurMoz96}
M.~Burger and Sh. Mozes, \emph{{{\rm CAT(-1)}}-spaces, divergence groups and
  their commensurators}, Journal of the American Mathematical Society
  \textbf{9} (1996), 57--93.

\bibitem[BM97]{BurMozCras}
M.~Burger and Sh. Mozes, \emph{Finitely presented simple groups and products of
  trees}, Comptes-Rendus de l'Acad\'emie des sciences de Paris \textbf{324}
  (1997), 747--752.

\bibitem[BM00a]{BurMozTrees}
M.~Burger and Sh. Mozes, \emph{Groups acting on trees: from local to global
  structure}, Publications Math\'ematiques de l'IH\'ES \textbf{92} (2000),
  113--150.

\bibitem[BM00b]{BurMozProd}
M.~Burger and Sh. Mozes, \emph{Lattices in product of trees}, Publications
  Math\'ematiques de l'IH\'ES \textbf{92} (2000), 151--194.

\bibitem[BM02]{BurMon}
M.~Burger and N.~Monod, \emph{Continuous bounded cohomology and applications to
  rigidity theory}, Geometric and Functional Analysis \textbf{12} (2002),
  no.~2, 219--280.

\bibitem[Bon03]{Bonvin}
P.~Bonvin, \emph{Strong boundaries and commensurator super-rigidity}, appendix
  to \cite{RemLin}, Institut Fourier preprint {\bf 590}, 2003.

\bibitem[Bou81]{Bou81}
N.~Bourbaki, \emph{Groupes et alg\`ebres de {L}ie {IV-VI}}, \'El\'ements de
  math\'ematique, Masson, 1981.

\bibitem[Bou97]{BouMostow}
M.~Bourdon, \emph{Immeubles hyperboliques, dimension conforme et rigidit\'e de
  {M}ostow}, Geometric and Functional Analysis \textbf{7} (1997), 245--268.

\bibitem[BP00]{BouPaj}
M.~Bourdon and Herv\'e Pajot, \emph{Rigidity of quasi-isometries for some
  hyperbolic buildings}, Commentarii Mathematici Helvetici \textbf{75} (2000),
  701--736.

\bibitem[Bro87]{BroCohomology}
K.S. Brown, \emph{Cohomology of groups}, Graduate texts in mathematics,
  vol.~82, Springer, 1987.

\bibitem[Bro89]{BroBuildings}
K.S. Brown, \emph{Buildings}, Springer, 1989.

\bibitem[BT72]{BruhatTits1}
F.~Bruhat and J.~Tits, \emph{Groupes r\'eductifs sur un corps local {I}.
  {D}onn\'ees radicielles valu\'ees}, Publications Math\'ematiques de l'IH\'ES
  \textbf{41} (1972), 5--251.

\bibitem[BT84]{BruhatTits2}
F.~Bruhat and J.~Tits, \emph{Groupes r\'eductifs sur un corps local {II}.
  {S}ch\'emas en groupes. {E}xistence d'une donn\'ee radicielle valu\'ee},
  Publications Math\'ematiques de l'IH\'ES \textbf{60} (1984), 5--184.

\bibitem[Bur95]{Burger}
M.~Burger, \emph{Rigidity properties of group actions on {{\rm
  CAT(0)}}-spaces}, Proceedings of the International Congress of Mathematicians
  (Z\"urich, August 1994) (S.D. Chatterji, ed.), Birkh\"auser Verlag, Basel,
  1995, pp.~761--769.

\bibitem[Cap03]{Caprace}
P.-E. Caprace, \emph{Isomorphisms of {K}ac-{M}oody groups}, M\'emoire,
  Universit\'e libre de Bruxelles, 2003.

\bibitem[Car72]{Carter}
R.W. Carter, \emph{Simple groups of lie type}, Wiley Intersciences, 1972.

\bibitem[Dav97]{Davis}
M.~Davis, \emph{Buildings are {{\rm CAT(0)}}}, Geometry and cohomology in group
  theory (P.H. Kropholler, G.A. Niblo, and R.~St\"ohr, eds.), LMS Lecture Notes
  Series 252, London Mathematical Society, Cambridge University Press, 1997,
  pp.~108--123.

\bibitem[DJ02]{DJ02}
J.~Dymara and T.~Januszkiewicz, \emph{Equivariant cohomology of buildings and
  of their automorphism groups}, Inventiones Mathematic\ae\, \textbf{150}
  (2002), no.~3, 579--627.

\bibitem[dlHV89]{HarVal}
P.~de~la Harpe and A.~Valette, \emph{La propri\'et\'e {{\rm (T)}} de {K}azhdan
  pour les groupes localement compacts}, Ast\'erisque, vol. 175, Soci\'et\'e
  Math\'ematique de France, 1989.

\bibitem[GdlH90]{GhyHar}
\'E. Ghys and P.~de~la Harpe, \emph{Sur les groupes hyperboliques d'apr\`es
  {M}ikhael {G}romov}, Progress in Mathematics, vol.~83, Birkh\"auser, 1990.

\bibitem[GJT98]{GJT}
Y.~Guivarc'h, L.~Ji, and J.C. Taylor, \emph{Compactifications of symmetric
  spaces}, Progress in Mathematics, vol. 156, Birkh\"auser, 1998.

\bibitem[Hag03]{HagComm}
F.~Haglund, \emph{Linearity and commensurability of uniform lattices of
  right-angled hyperbolic buildings}, private communication, 2003.

\bibitem[Har69]{Harder}
G.~Harder, \emph{Minkowskische {R}eduktionstheorie \"uber
  {F}unktionenk\"orpern}, Inventiones Mathematic\ae\, \textbf{7} (1969),
  33--54.

\bibitem[HP98]{HagPau}
F.~Haglund and F.~Paulin, \emph{Simplicit\'e de groupes d'automorphismes
  d'espaces \`a courbure n\'egative}, Geom. Topol. Monogr. \textbf{1} (1998),
  181--248.

\bibitem[Jac89]{Jacobson}
N.~Jacobson, \emph{Basic {A}lgebra {II} (2nd edition)}, Freeman and Co, 1989.

\bibitem[Kai02]{Kaim}
V.A. Kaimanovich, \emph{Double ergodicity of the {P}oisson boundary and
  applications to bounded cohomology}, to appear in Geometric and Functional
  Analysis, 2002.

\bibitem[KP85]{KacPet}
V.~Kac and D.~Peterson, \emph{Defining relations for certain
  infinite-dimensional groups}, \'Elie Cartan et les math\'ematiques
  d'aujourd'hui. The mathematical heritage of \'Elie Cartan. Lyon, 25-29 juin
  1984 (Soci\'et\'e math\'ematique~de France, ed.), Ast\'erisque Hors-S\'erie,
  1985, pp.~165--208.

\bibitem[Liu94]{Liu}
Y.~Liu, \emph{Density of the commensurability group of uniform tree lattices},
  Journal of Algebra \textbf{165} (1994), 346--359.

\bibitem[LMZ94]{LMZ94}
A.~Lubotzky, Sh. Mozes, and R.J. Zimmer, \emph{Superrigidity for the
  commensurability group of tree lattices}, Commentarii Mathematici Helvetici
  \textbf{69} (1994), 523--548.

\bibitem[Mar91]{Margulis}
G.A. Margulis, \emph{{D}iscrete {S}ubgroups of {S}emisimple {L}ie {G}roups},
  Ergebnisse der Mathematik und ihrer Grenzgebiete (3), vol.~17, Springer,
  1991.

\bibitem[Mar94]{Marg94}
G.A. Margulis, \emph{Superrigidity of commensurability subgroups and
  generalized harmonic maps}, unpublished, 1994.

\bibitem[Mas88]{Maskit}
B.~Maskit, \emph{Kleinian groups}, Grundlehren der matematischen
  Wissenschaften, vol. 287, Springer, 1988.

\bibitem[Mou88]{Moussong}
G.~Moussong, \emph{Hyperbolic {C}oxeter groups}, PhD thesis, Ohio State
  University, 1988.

\bibitem[Moz99]{MozNagao}
Sh. Mozes, \emph{Trees, lattices and commensurators}, Algebra, $K$-Theory,
  Groups, and Education: On the Occasion of Hyman Bass's 65th Birthday (T.Y.
  Lam and A.R. Magid, eds.), Contemporary Mathematics, no. 243, American
  Mathematical Society, 1999, pp.~145--151.

\bibitem[MR95]{MuhRon}
B.~M\"uhlherr and M.A. Ronan, \emph{Local to global structure in twin
  buildings}, Inventiones Mathematic\ae\, \textbf{122} (1995), 71--81.

\bibitem[Mue99]{Muhlherr}
B.~Muehlherr, \emph{Locally split and locally finite buildings of $2$-spherical
  type}, J. Reine angew. Math. \textbf{511} (1999), 119--143.

\bibitem[PR94]{PR94}
V.~Platonov and A.~Rapinchuk, \emph{Algebraic groups and number theory},
  Academic press, 1994.

\bibitem[Pra77]{Prasad}
G.~Prasad, \emph{Strong approximation for semi-simple groups over function
  fields}, Annals of Mathematics \textbf{105} (1977), 553--572.

\bibitem[Rem99]{RemCras}
B.~Remy, \emph{Construction de r\'eseaux en th\'eorie de {K}ac-{M}oody},
  Comptes-Rendus de l'Acad\'emie des sciences de Paris \textbf{329} (1999),
  475--478.

\bibitem[Rem02a]{RemNewton}
B.~Remy, \emph{Classical and non-linearity properties of {K}ac-{M}oody groups},
  Rigidity in Dynamics and Geometry (Newton Institute, 2000) (M.~Burger and
  A.~Iozzi, eds.), Springer Verlag, 2002, pp.~391--406.

\bibitem[Rem02b]{RemAst}
B.~Remy, \emph{Groupes de {K}ac-{M}oody d\'eploy\'es et presque d\'eploy\'es},
  Ast\'erisque, vol. 277, Soci\'et\'e Math\'ematique de France, 2002.

\bibitem[Rem02c]{RemGD}
B.~Remy, \emph{Immeubles de {K}ac-{M}oody hyperboliques, groupes non isomorphes
  de m\^eme immeuble}, Geometri\ae \ Dedicata \textbf{90} (2002), 29--44.

\bibitem[Rem03a]{RemBie}
B.~Remy, \emph{{K}ac-{M}oody groups: split and relative theories. {L}attices},
  Groups: Topological, Combinatorial and Arithmetic Aspects (Bielefeld, 1999)
  (T.W. M\"uller, ed.), LMS Lecture Notes Series, London Mathematical Society,
  Cambridge University Press, 2003, pp.~401--445.

\bibitem[Rem03b]{RemImage}
B.~Remy, \emph{Linear images of {K}ac-{M}oody groups}, in preparation, 2003.

\bibitem[Rem03c]{RemLin}
B.~Remy, \emph{Topological simplicity, commensurator super-rigidity and
  non-linearities of {K}ac-{M}oody groups}, Institut Fourier preprint {\bf
  590}, 2003.

\bibitem[Ron86]{RonFree}
M.A. Ronan, \emph{A construction of buildings with no rank 3 residues of
  spherical type}, Buildings and the geometry of diagrams (Lectures given at
  the third 1984 session of the Centro Internazionale Matematico Estivo (CIME)
  held in Como, August 26-September 4, 1984) (L.A. Rosati, ed.), Lecture Notes
  in Mathematics, vol. 1181, Springer Verlag, Berlin, 1986, pp.~242--248.

\bibitem[Ron89]{RonLec}
M.A. Ronan, \emph{Lectures on {B}uildings}, Perspectives in Mathematics,
  vol.~7, Academic Press, 1989.

\bibitem[RR02]{RemRon}
B.~Remy and M.~Ronan, \emph{Topological groups of {K}ac-{M}oody type,
  {F}uchsian twinnings and their lattices}, Institut Fourier preprint {\bf
  563}, 2002.

\bibitem[Sha00]{Shalom}
Y.~Shalom, \emph{Rigidity of commensurators and irreducible lattices},
  Inventiones Mathematic\ae\, \textbf{141} (2000), 1--54.

\bibitem[Tit70]{TitsTree}
J.~Tits, \emph{Sur le groupe des automorphismes d'un arbre}, Essays on topology
  and related topics (M\'emoires d\'edi\'es \`a Georges de Rham (A.~H\ae fliger
  and R.~Narasimhan, eds.), Springer Verlag, 1970, pp.~188--211.

\bibitem[Tit72]{Tit72}
J.~Tits, \emph{Free subgroups in linear groups}, Journal of Algebra \textbf{20}
  (1972), 250--270.

\bibitem[Tit74]{TitsSpherical}
J.~Tits, \emph{Buildings of spherical type and finite {BN}-pairs}, Lecture
  Notes in Mathematics, vol. 386, Springer, 1974.

\bibitem[Tit87]{Tit87}
J.~Tits, \emph{Uniqueness and presentation of {K}ac-{M}oody groups over
  fields}, Journal of Algebra \textbf{105} (1987), 542--573.

\bibitem[Tit89]{Tit89}
J.~Tits, \emph{Groupes associ\'es aux alg\`ebres de {K}ac-{M}oody}, S\'eminaire
  Bourbaki 700, Ast\'erisque, vol. 177-178, Soci\'et\'e Math\'ematique de
  France, 1989, pp.~7--31.

\bibitem[Tit90]{TitsCDF}
J.~Tits, \emph{Th\'eorie des groupes}, Annuaire Coll\`ege de France (R\'esum\'e
  de cours), 1990, pp.~87--103.

\bibitem[Tit92]{Tit92}
J.~Tits, \emph{Twin buildings and groups of {K}ac-{M}oody type}, Groups,
  Combinatorics \& Geometry (LMS Symposium on Groups and Combinatorics, Durham,
  July 1990) (M.~Liebeck and J.~Saxl, eds.), LMS Lecture Notes Series, vol.
  165, Cambridge University Press, 1992, pp.~249--286.

\bibitem[Ven88]{Venka}
T.N. Venkataramana, \emph{On superrigidity and arithmeticity of lattices in
  semisimple groups over local fields of arbitrary characteristic}, Inventiones
  Mathematic\ae\, \textbf{92} (1988), no.~2, 255--302.

\bibitem[Zim84]{Zimmer}
R.J. Zimmer, \emph{Ergodic {T}heory and {S}emisimple {G}roups}, Monographs in
  Mathematics, vol.~81, Birkh\"auser, 1984.

\end{thebibliography}

\vspace{1cm}

\addtolength{\parindent}{-1.6pt} 

\gec 

Institut Fourier -- UMR CNRS 5582\\
Universit\'e de Grenoble 1 -- Joseph Fourier\\
100, rue des maths -- BP 74\\
38402 Saint-Martin-d'H\`eres Cedex -- France\\
{\tt bremy@fourier.ujf-grenoble.fr}

\end{document}